\documentclass{amsart} 
\usepackage{amssymb,latexsym,amsfonts,amsmath}
\usepackage{url}
\usepackage{epsfig}

\newcommand{\R}{{\mathbb{R}}}

\newcommand{\N}{{\mathbb{N}}}

\newcommand{\ie}{{\it i.e. }}
\newcommand{\eg}{{\it e.g. }}

\usepackage{color} 

\newtheorem{thm}{Theorem}
\newtheorem{lem}{Lemma}

\newtheorem{prob}{Problem}
\newtheorem{assum}{Assumption}
\newtheorem{prop}{Proposition}
\newtheorem{cor}{Corollary}

\newtheorem{defn}{Definition}

\newtheorem{rem}{Remark}
\newenvironment{pf}{\begin{proof}}{\end{proof}}

\begin{document}

 \title{Asynchronous Decentralized Event-triggered Control}
\author{Manuel Mazo Jr. and Ming Cao}
 \thanks{This work was partially supported by grants from the Dutch Organization for Scientific Research (NWO), the Dutch Technology Foundation (STW) and the European Union Seventh Framework Programme [FP7/2007-2013]  under grant agreement 257462 HYCON2 Network of Excellence.}
 \thanks{M. Mazo Jr is with INCAS$^3$, Assen and the Department of Discrete Technology and Production Automation, University of Groningen, The Netherlands, 
{\tt\small m.mazo@rug.nl}}
 \thanks{M. Cao is with the Department of Discrete Technology and Production Automation, University of Groningen, The Netherlands, 
{\tt\small m.cao@rug.nl}
}

\maketitle

\begin{abstract}                          
In this paper we propose an approach to the implementation of controllers with decentralized strategies triggering controller updates. We consider set-ups with a central node in charge of the computation of the control commands, and a set of not co-located sensors providing measurements to the controller node. The solution we propose does not require measurements from the sensors to be synchronized in time. The sensors in our proposal provide measurements in an aperiodic way triggered by local conditions. Furthermore, in the proposed implementation (most of) the communication between nodes requires only the exchange of one bit of information (per controller update), which could aid in reducing transmission delays and as a secondary effect result in fewer transmissions being triggered.
\end{abstract}

\section{Introduction}

Aperiodic control techniques have recently gained much attention due to the opportunities they open to reduce bandwidth and computation requirements in control systems implementations~\cite{Astrom,Tabuada,AntaTAC}. These savings are especially relevant in the implementation of control loops over wireless channels~\cite{AraujoDCOSS,Rabi08}. In those set-ups there is not only a limited bandwidth available, but also sensor nodes may have limited energy provided by batteries. It is therefore interesting to explore approaches which may save energy expenditures at the sensors by \eg reducing the number of transmissions from those sensors, or reducing the amount of time the sensor nodes need to keep their radios listening for possible communications from other nodes. While there is an extensive recent literature on event-triggered control aimed at reducing the amount of transmissions necessary to close the control loop while maintaining stability~\cite{Heemels08,Cervin08,Lunze10,Molin10,Wang11}, the problem of reducing listening time has received less attention~\cite{WeimerADHS12,Donkers12,MazoCaoADHS2012,MazoCaoCDC2011}. Nevertheless, it is a well-known phenomena in the sensor networks community that reducing listening times has a bigger impact on the power burden than reducing transmissions~\cite{Estrin02}. In the present paper we try to bridge this gap by proposing controller implementations focused on reducing listening times. In order to attain this goal, we propose a technique in which the sensors do not need to coordinate with each other, and therefore do not need to listen to each other. Instead, in the proposed implementation the sensors send measurements triggered by local conditions, irrespective of what the other sensors are doing, in contrast with previous work on decentralized triggering~\cite{MazoTACWSAN}. With respect to other work on decentralized or distributed event-triggered control, we do not impose any weak coupling assumptions or very restrictive dynamics, as is often the case in work focused on multi-agent systems~\cite{Wang11,Dimarogonas12}. Note that also~\cite{Wang11,Dimarogonas12} suffer from the drawback of continuous listening. Arguably, the work closest to ours is that presented in~\cite{Donkers12}, however restricted to the study of linear systems.

The implementation that we propose also enables the stabilization of systems employing communication packets with very small payload. In particular, our technique reduces the amount of payload needed to essentially one bit. To appreciate the relevance of reducing the packets payload, besides reducing power consumption~\cite{Estrin02}, one must notice that a large portion of delays in communications are due to transmission delay. These transmission delays are dependent on the size of the packages transmitted, and thus reducing the payload will indirectly reduce the communication delays present in the system. Event-triggered implementations of control systems accommodate delays by making more conservative the conditions that trigger communications than those in the delay free case. Employing more conservative conditions results, in general, in more frequent transmissions of measurements. Thus, a reduction on the payload is also expected to result in a reduction on the amount of transmissions from the sensors to the controller.

The ideas in the present paper will remind the reader of dynamic quantizers for control~\cite{DynamicQuantizers} and of dead-band control~\cite{DeadbandTilbury}. We have, in a way, combined those ideas with recent approaches to event-triggered control stemming from~\cite{Tabuada} to provide a formal analysis of implementations benefiting from all those ideas. The current paper is the result of merging previous conference contributions by the authors, providing a unified analysis and removing early mistakes and imprecise statements. As such it should be seen as a more accurate and easier to follow analysis of the proposals in~\cite{MazoCaoCDC2011} and~\cite{MazoCaoADHS2012}.

The remainder of this paper is organized as follows: Section~\ref{sec:prelim} introduces some mathematical preliminaries and notation; this is followed by a formalization of the problem and an introduction to the idea of asynchronous event-triggered control in Section~\ref{sec:review_adet}; Section~\ref{sec:contrib} contains our main contribution, in the form of an asynchronous event-triggered implementation attaining asymptotic stability with one bit communication exchange; the paper concludes with two examples in Section~\ref{sec:example} and a discussion in Section~\ref{sec:discuss}.

\section{Preliminaries}
\label{sec:prelim}

We denote the positive
real numbers by $\R^+$ and by
$\R_0^+=\R^+\cup\lbrace0\rbrace$. We also use $\N$ to denote the natural numbers including zero. 
The usual Euclidean ($l_2$) vector norm is represented by $|\cdot|$. 
 When applied to a matrix $|\cdot|$ denotes the
 $l_2$ induced matrix norm. A matrix
 $P\in\R^{n\times n}$ is said to be positive definite, denoted by $P>0$, whenever
 $x^TPx>0$ for all $x\neq 0$, $x\in\R^n$. By $\lambda_m(P),\lambda_M(P)$
 we denote the minimum and maximum eigenvalues of $P$ respectively. 
A function $f:\R^n\to\R^m$ is said to be Lipschitz
continuous on compacts if for every compact set $S\subset\R^n$ there exists a constant $L\in\R_0^+$ such that:
$|f(x)-f(y)|\leq L|x-y|, \; \forall x,y \in S$. 
For a function $f:\R^n\to\R^n$ we denote by $f_i:\R^n\to\R$ the function whose image is the projection of $f$ on its i-th coordinate, \ie $f_i(x)=\Pi_i(f(x))$. Consequently, given a Lipschitz continuous function $f$, we also denote by $L_{f_i}$ the Lipschitz constant of $f_i$. 
A function $\gamma:\R^+_0\to\R^+_0$, is of class~$\mathcal{K}_\infty$ if it is
continuous, strictly increasing, $\gamma(0)=0$ and $\gamma(s)\to\infty$ as
$s\to\infty$. 
A continuous function $\beta:\R^+_0\times\R^+_0\to\R^+_0$ is of
class~$\mathcal{KL}$ if $\beta(\cdot,\tau)$ is of class~$\mathcal{K}_\infty$
for each $\tau\geq0$ and $\beta(s,\cdot)$ is monotonically decreasing to zero
for each $s\geq0$. 
Given an essentially bounded function $\delta:\R^+_0\to\R^m$ we denote by $\Vert
\delta \Vert_\infty$ the $\mathcal{L}_\infty$ norm, \ie $\Vert \delta
\Vert_\infty=\text{ess}\sup_{t\in\R^+_0}\lbrace{|\delta(t)|\rbrace}$.

We discuss the following properties of a closed-loop system in this paper:

\begin{defn}[Asymptotic Stability]$\,$\\
A system $\dot{\xi}(t)=f(\xi(t)),\; t\in\R_0^+,\, \xi(t)\in\R^n$ is said to be uniformly globally asymptotically stable (UGAS) if there exists $\beta\in\mathcal{KL}$ such that for any $t_0\geq 0$ the following holds:
\begin{eqnarray*}
\forall \xi(t_0)\in\R^n\,\Rightarrow |\xi(t)|\leq \beta(|\xi(t_0)|, t-t_0),\;\forall t\geq t_0;
\end{eqnarray*}
\end{defn}
\begin{defn}[Practical Stability]$\,$\\
A system $\dot{\xi}(t)=f(\xi(t)),\; t\in\R_0^+,\, \xi(t)\in\R^n$ is said to be  uniformly globally practically stable (UGPS) if there exist $\beta\in\mathcal{KL}$ and $\delta>0$ such that for any $t_0\geq 0$ the following holds:
\begin{eqnarray*}
\forall \xi(t_0)\in\R^n\,\Rightarrow |\xi(t)|\leq \beta(|\xi(t_0)|, t-t_0)+\delta,\;\forall t\geq t_0;
\end{eqnarray*}
\end{defn}

We only consider asymptotic stability with respect to the origin, and not with respect to sets. On the other hand, our notion of practical stability can be identified with asymptotic stability to a set: a ball around the origin.
Furthermore, we only consider global properties for simplicity of the presentation. Local versions of the results should not be much harder to provide following the same reasoning we follow in this paper.

The notion of Input-to-State stability (ISS)~\cite{Sontag:2008p2453} 
will be central to our discussion: 
\begin{defn}[Input-to-State Stability]
\quad\\
A control system $\dot{\xi}=f(\xi,\upsilon)$ is said to be (uniformly globally) Input-to-state stable (ISS) with respect to $\upsilon$ if there exists $\beta\in\mathcal{KL}$, $\gamma\in\mathcal{K}_\infty$ such that for any $t_0\in\R^+_0$ the following holds:
\begin{eqnarray*}
\forall \xi(t_0)\in\R^n,\, \Vert \upsilon \Vert_\infty<\infty \Rightarrow\\
|\xi(t)|\leq\beta(|\xi(t_0)|,t-t_0)+\gamma(\Vert
\upsilon \Vert_\infty),\,\forall t\geq t_0.
\end{eqnarray*}
\end{defn}

Rather than using its definition, in our arguments we rely on the following characterization: a system is ISS if and only if there exists an ISS Lyapunov function~\cite{Sontag:2008p2453}.
\begin{defn}[ISS Lyapunov function]
\quad\\
A smooth function $V:\R^n\to\R_0^+$ is said to be an ISS Lyapunov function for the closed-loop system $\dot\xi=f(\xi,\upsilon)$ if there exists class $\mathcal{K}_\infty$ functions $\underline{\alpha}$, $\overline{\alpha}$, $\alpha_v$ and $\alpha_e$ such that for all $x\in\R^n$ and $u\in\R^m$ the following is satisfied:
\begin{eqnarray}
\label{eq:V_bounded}\notag
\underline{\alpha}(|x|)\leq& V(x)&\leq \overline{\alpha}(|x|)\\
\nabla V\cdot f(x,u)&\leq& -\alpha_v\circ V(x)+\alpha_e(|u|)
\end{eqnarray}
\end{defn}

Often we use the shorthand $\dot{V}(x,u)$ to denote the Lie derivative $\nabla V\cdot f(x,u)$, and $\circ$ to denote function composition, \ie $f\circ g(t)=f(g(t))$.

Finally, we employ the following, rather trivial, result in some of our arguments:

\begin{lem}
\label{lem:growth_bound}
Given two $\mathcal{K}_\infty$ functions $\alpha_1$ and $\alpha_2$, there exists some constant $L<\infty$ such that:
$$
\lim_{s\to 0}\frac{\alpha_1(s)}{\alpha_2(s)}=L
$$
if and only if  for all $S<\infty$ there exists a positive $\kappa<\infty$ such that:
$$
\forall s\in ]0,S],\, \alpha_1(s)\leq \kappa \alpha_2(s).
$$
\end{lem}
\begin{pf}
The necessity side of the equivalence is trivial, thus we concentrate on the sufficiency part. By assumption, we know that the limit of the ratio of the functions tends to $L$ as $s\to 0$, and therefore, $\forall \epsilon>0$, $\exists \delta>0$ such that $\alpha_1(s)/\alpha_2(s)<L+\epsilon$ for all $s\in]0,\delta[$. As $\alpha_1,\alpha_2\in\mathcal{K}_\infty$ we know that in any compact set excluding the origin the function $\alpha_1(s)/\alpha_2(s)$ is continuous and therefore attains a maximum, \ie $\alpha_1(s)/\alpha_2(s)<M,\, \forall s\in[\delta,S],\, 0<\delta<S$. Putting these two results together we have that $
\forall s\in ]0,S], \, S<\infty,\; \alpha_1(s)\leq \kappa \alpha_2(s)$, where $\kappa=\max\lbrace L+\epsilon, M \rbrace$.
\end{pf}

\section{Asynchronous Event-triggered Control}
\label{sec:review_adet}

The problem we aim at solving is that of controlling systems of the form:
\begin{equation}
\label{eq:system}
\dot{\xi}(t)=f(\xi(t),\upsilon(t)),\qquad \forall t\in \R^+_0,
\end{equation}
where $\xi:\R^+_0\to\R^n$ and $\upsilon:\R^+_0\to\R^m$ and the full state is assumed to be measured.
In particular, we are interested in finding stabilizing sample-and-hold implementations of a controller $\upsilon(t)=k(\xi(t))$ such that
updates can be performed transmitting asynchronous and aperiodic measurements of each entry of the state vector. 
Furthermore, if possible we would like to do so while reducing the amount of transmissions.
This problem can be formalized as follows:
\begin{prob}
\label{prb:main}
Given system (\ref{eq:system}) and a controller $k:\R^n\to\R^m$
find sequences of update times $\lbrace{t^i_{r_i}\rbrace},\, r_i\in\N$ for each sensor $i=1,\ldots,n$
such that an asynchronous sample-and-hold controller implementation:
\begin{eqnarray}
\label{eq:controller}
\upsilon_j(t)=&k_j(\xi_1(t^1_{r_1}), \xi_2(t^2_{r_2}),\ldots, \xi_n(t^n_{r_n})),\; j=1,\ldots,m\\\notag 
& t\in[\max_{i=1,\ldots,n}\lbrace t^i_{r_i}\rbrace,\min_{i=1,\ldots,n}\lbrace t^i_{r_i+1}\rbrace[,
\end{eqnarray}
renders the closed-loop system UGAS (or UGPS).
\end{prob}

Let us start introducing the two main technical assumptions of the paper:
\begin{assum}[ISS w.r.t. measurement errors]
\label{ass:ISS}
Given system~(\ref{eq:system}) there exists a controller $k:\R^n\to\R^m$ such that the closed-loop
system
\begin{equation}
\label{eq:sys_iss}
\dot{\xi}(t)=f(\xi(t),k(\xi(t)+\varepsilon(t)),\qquad \forall t\in \R^+_0
\end{equation}
is ISS with respect to measurement errors $\varepsilon$.
\end{assum}

\begin{assum}[Lipschitz continuity]
\label{ass:lipschitz}
\quad\\
The functions $f$ and $k$, defining the dynamics and controller of a system, are locally Lipschitz on compacts.
\end{assum}
Note that this last assumption guarantees the existence and uniqueness of solutions of the closed-loop system.

Representing the effect of the sample-and-hold as a measurement error at each sensor for all $i=1,\ldots,n$ as:
$$
\varepsilon_i(t)=\xi_i(t^i_{r_i})-\xi_i(t),\; t\in[t^i_{r_i}, t^i_{r_i+1}[,\;r_i\in\N
$$
we propose rules defining implicitly the sequences of update times $\lbrace{t^i_{r_i}\rbrace}$ for each sensor $i$:
\begin{equation}
\label{eq:upd_rule}
t^i_{r_i}:=\min \lbrace t>t^i_{r_i-1}\, |\, \varepsilon^2_i(t)=\eta_i  \rbrace,
\end{equation}
where $\eta_i>0$ are design parameters.

For convenience and compactness of the presentation we introduce the new variable: 
\begin{small}
\begin{equation}
\label{eq:relation_eta}
\eta=\sqrt{\sum_{i=1}^n \eta_i},
\end{equation}
\end{small}
and consider it as a design parameter that once specified restricts the choices of $\eta_i$ to be used at each sensor.
We remark now that with this definition the update rule~(\ref{eq:upd_rule}) implies that $|\varepsilon(t)|\leq\eta$ (with equality attained only when all sensors trigger simultaneously).
Furthermore, we assume that the local parameters $\eta_i$ are defined through an appropriate scaling: 
\begin{equation}
\label{eq:eta_i}
\eta_i:=\theta_i^2\eta^2,\, |\theta|=1,
\end{equation} 
with $\theta_i$ as design constants.

Now, we can state the following lemma which will be central in the rest of the discussion:

\begin{lem}[Inter-transmission times bound]
\label{lem:bounds}
\quad\\
If Assumptions~\ref{ass:ISS} and \ref{ass:lipschitz} hold,
for any $\eta>0$,  a lower bound for the minimum time between transmissions of a sensor, for all time $t\geq t_0$, is given by:
\begin{equation}
\label{eq:min_time_constant}
\tau_i^*:= L_{f_i}^{-1}\theta_i\frac{\eta}{\eta+\underline{\alpha}^{-1}(\max\lbrace V(\xi(t_0)),\alpha_v^{-1}\circ\alpha_e(\eta) \rbrace)}.
\end{equation}
where $L_{f_i}$ denotes the Lipschitz constant of the function $f_i(x,k(x+e))$ for $|x|\leq\underline{\alpha}^{-1}(\max\lbrace V(\xi(t_0)),\alpha_v^{-1}\circ\alpha_e(\eta) \rbrace)$ and $|e|\leq \eta$.
\end{lem}
\begin{pf}
Let us denote in what follows by:\newline $\overline{f}_i(y,z)=\max_{(x,e)\in S(y,z)}|f_i(x,k(x+e))|$ and by $S(y,z)=\lbrace (x,e)\in\R^{n\times n}\,|\, V(x)\leq y,\,|e|\leq z\rbrace$.
Start by recalling that the minimum time between events at a sensor is given by the time it takes for $|\varepsilon_i|$ to evolve from the value $|\varepsilon_i(t^i_{k_i})|=0$
to $|\varepsilon_i(t^{i-}_{k_i+1})|=\sqrt{\eta_i}$. Therefore all that needs to be proved is the existence of an upper bound on the rate of change of $|\varepsilon_i|$.
By Assumption~\ref{ass:ISS} we have that: $|e|\leq\eta,\, V(x)\geq\alpha_v^{-1}\circ\alpha_e(\eta)\Rightarrow \dot{V}(x,e)\leq 0$ and thus
we have that for any $r\in\R^+_0$ the set $S(\alpha_v^{-1}\circ\alpha_e(\eta)+r,\eta)$ is forward invariant.
One can trivially bound the evolution of $|\varepsilon|$ as:
\begin{eqnarray*}
\frac{d}{dt}|\varepsilon_i|&\leq&|\dot{\varepsilon}_i|=|f_i(\xi,k(\xi+\varepsilon))|,
\end{eqnarray*}
and the maximum rate of change of $|\varepsilon_i|$ by\newline $\overline{f}_i(\max\lbrace V(\xi(t_0)),\alpha_v^{-1}\circ\alpha_e(\eta) \rbrace,\eta).$ 
Note that the existence of such a maximum is guaranteed by the continuity of the maps $f$ and $k$ and the compactness of the set $S(\max\lbrace V(\xi(t_0)),\alpha_v^{-1}\circ\alpha_e(\eta) \rbrace,\eta)$.
Assumption~\ref{ass:lipschitz} implies that $f_i(x,k(x+e))$ is also locally Lipschitz, and thus one can further bound
$\overline{f}_i(\max\lbrace V(\xi(t_0)),\alpha_v^{-1}\circ\alpha_e(\eta) \rbrace,\eta) \leq$ $L_{f_i}\left(\underline{\alpha}^{-1}(\max\lbrace V(\xi(t_0)),\alpha_v^{-1}\circ\alpha_e(\eta) \rbrace)+\eta\right).$ 
Finally, recalling that $\eta_i=\theta_i^2\eta^2$, a lower bound for the inter-transmission times is given by~(\ref{eq:min_time_constant})
which proves the statement.
\end{pf}

With this result one can state the following theorem:

\begin{thm}[UGPS]
\label{thm:pract}
\quad\\
If Assumptions~\ref{ass:ISS}~and~\ref{ass:lipschitz} hold,
then the closed-loop system~(\ref{eq:system}),~(\ref{eq:controller}),~(\ref{eq:upd_rule}) is UGPS.
Moreover, the time between transmissions of measurements at each sensor is bounded from below by some $\tau^*_i>0$.
\end{thm}
\begin{pf}
From Lemma~\ref{lem:bounds} we know that there exists some minimum time between the triggering of events at the different sensors. This, together with that the number of sensors is finite, guarantees that there are no Zeno executions of the closed-loop system. Assumption~\ref{ass:ISS} provides the bound: $|\xi(t)|\leq\beta(|\xi(t_0)|,t-t_0)+\gamma(\Vert \varepsilon \Vert_\infty)$, and the proposed implementation forces $\Vert \varepsilon \Vert_\infty\leq\eta$. Thus, using these two bounds, and ruling out any possible Zeno solution, let us state that $|\xi(t)|\leq\beta(|\xi(t_0)|,t-t_0)+\delta$, where $\delta:=\gamma(\eta)$, which finalizes the proof.
\end{pf}

In general, employing a constant threshold value $\eta$ establishes a trade-off between the size of the inter-transmission times and the size of the set to which the system converges. In order to solve this problem the following lemma establishes an approach to construct asymptotic stabilizing asynchronous implementations. The main idea is to let the parameter $\eta$ change over time as:
\begin{equation}
\label{eq:eta_variable}
\eta(t)=\eta(t^c_{r_c}),\;\,\forall t\in[t^c_{r_c},t^c_{r_c+1}[,
\end{equation}
with $\lbrace t^c_{r_c}\rbrace$ a divergent sequence of times (to be defined later)
and use the local update rules:
\begin{eqnarray}
\label{eq:upd_rule_dyn}\notag
t^i_{r_i}&:=&\min \lbrace t>t^i_{r_i-1}\, |\, \varepsilon^2_i(t)=\eta_i(t)\rbrace,\\
\eta_i(t)&:=&\theta_i^2\eta(t)^2,\, |\theta|=1,
\end{eqnarray}
\begin{lem}[UGAS]
\label{lem:UGAS}
\quad\\
If Assumptions~\ref{ass:ISS}~and~\ref{ass:lipschitz} are satisfied and the following two conditions hold:
\begin{itemize}
 \item $\lim_{r_c\to\infty}\eta(t^c_{r_c})\to 0$;
 \item $\frac{\underline{\alpha}^{-1}(\max\lbrace V(\xi(t^c_{r_c})),\alpha_v^{-1}\circ\alpha_e(\eta(t^c_{r_c})) \rbrace)}{\eta(t^c_{r_c})}\leq\kappa<\infty$ for all $t^c_{r_c}$.
\end{itemize}
for a divergent sequence of times $\lbrace t^c_{r_c}\rbrace$, the closed loop system~(\ref{eq:system}), (\ref{eq:controller}), (\ref{eq:eta_variable}), (\ref{eq:upd_rule_dyn})  is UGAS.

\end{lem}
\begin{pf}
In view of Lemma~\ref{lem:bounds} the second condition of this lemma guarantees that there exists a minimum time between events at each sensor for every time interval $[t^c_{r_c}, t^c_{r_c+1}[$, \ie $t^i_{r_i+1}-t^i_{r_i}>\tau_i^*$ for all $i=1,\ldots,n$ and $t^i_{r_i}, t^i_{r_i+1}\in[t^c_{r_c}, t^c_{r_c+1}[$. It could happen however that some sensor update coincides with an update of the thresholds, \ie $t^i_{r_i+1}=t^c_{r_c+1}$, which could lead to two arbitrarily close events of sensor $i$. Similarly, events from two different sensors could be generated arbitrarily close to each other. Nonetheless, as the sequence $\lbrace t^c_{r_c}\rbrace$ is divergent (by assumption), and there is a finite number of sensors, none of these two effects can lead to Zeno executions. Thus the second assumption guarantees the exclusion of any possible Zeno behavior of the system.

Now note that from Assumption~\ref{ass:ISS} we have that for all $t\in[t^c_{r_c}, t^c_{r_c+1} [$ the following bound holds: $V(\xi(t))\leq \max\lbrace V(\xi(t^c_{r_c})), \alpha_v^{-1}\circ\alpha_e(\eta(t^c_{r_c}))\rbrace$. Furthermore, the second condition also implies that $V(\xi(t^c_{r_c}))\leq \gamma_1(\eta(t^c_{r_c}))$ for some $\gamma_1\in\mathcal{K}_\infty$, which can be shown constructively: the second condition implies that at $t^c_{r_c}$ either $V(\xi(t^c_{r_c}))\leq \alpha_v^{-1}\circ\alpha_e(\eta(t^c_{r_c}))$ or $V(\xi(t^c_{r_c}))\leq \underline{\alpha}(\kappa\eta(t^c_{r_c}))$, and thus $\gamma_1(s)=\max\lbrace\alpha_v^{-1}\circ\alpha_e(s), \underline{\alpha}(\kappa s) \rbrace$ provides the desired result. 

With these last two pieces of information, one can bound for all $t\in[t^c_{r_c}, t^c_{r_c+1}[$: 
$$V(\xi(t))\leq \max\lbrace \gamma_1(\eta(t^c_{r_c})), \alpha_v^{-1}\circ\alpha_e(\eta(t^c_{r_c}))\rbrace:=\gamma_2(\eta(t^c_{r_c})),$$
where $\gamma_2\in\mathcal{K}_\infty$. Next we notice that the first condition of the Lemma implies that $\exists \tilde{\beta}_\eta\in\mathcal{KL}$ such that ${\eta}(t)\leq \tilde{\beta}_\eta(\eta(t^c_0), t-t_0)$ for all $t\geq t^c_0$. Putting together these last two bounds, and assuming that $\kappa_0:=\eta(t^c_0)/V(\xi(t^c_0)<\infty$ one can conclude that 
$$V(\xi(t))\leq \tilde\beta(\gamma_2(\kappa_0V(t^c_{0})), t-t_0),\,\forall t\geq t^c_0.$$
Finally this last bound guarantees that
\begin{eqnarray}
|\xi(t)|&\leq& \underline\alpha^{-1}(\tilde\beta(\gamma_2(\kappa_0 \overline\alpha(|\xi(t^c_{0})|)), t-t_0))\\
&:=&\beta(|\xi(t^c_0)|, t-t_0),\, \forall t\geq t^c_0
\end{eqnarray}
with $\beta\in\mathcal{KL}$ which finalizes the proof.
\end{pf}

\section{An UGAS asynchronous implementation}
\label{sec:contrib}

In this section we present an implementation for asymptotic stability employing only asynchronous measurements from all sensors. This approach relies on predefining an update policy for the time-varying threshold $\eta(t^c_{r_c})$, active in the interval $t\in[t^c_{r_c}, t^c_{r_c+1}[$, as follows:
\begin{eqnarray}
\label{eq:eta_upd_rule}
\eta(t)&=&\eta(t^c_{r_c}),\, t\in[t^c_{r_c}, t^c_{r_c+1}[\\
\eta(t^c_{r_c+1})&=&\mu\eta(t^c_{r_c}),
\end{eqnarray}
for some $\mu\in]0,1[$. 
Given this update policy one can design an event-triggered policy to decide the sequence of times $\lbrace t^c_{r_c}\rbrace$ such that the system is rendered asymptotically stable. Furthermore, as we show later in this section, such a fully event-triggered implementation enables asymptotic implementations only requiring the exchange of one bit of information whenever communication between a sensor and controller, and vice-versa, is necessary.
This new strategy uses two independent triggering mechanisms:
\begin{itemize}
 \item {\bf Sensor-to-controller:} Sensors send measurements to the controller whenever the local threshold is violated. As explained in the previous section, the update of the control commands is done with the measurements as they arrive in an asynchronous fashion.
 \item {\bf Controller-to-sensor:} The controller commands the sensors to reduce the threshold used in their triggering condition when the system has "slowed down" enough to guarantee that the inter-sample times remain bounded from below. The triggering condition is based on the system entering a smaller Lyapunov level set which, by the Lipschitz assumption on the dynamics, implies slower dynamics of the system. The controller checks this condition only in a periodic fashion, with period $\tau^c$, and therefore the sensors only need to listen at those time instants.
\end{itemize}
The mechanism to trigger sensor to controller communication has already been analyzed in Section~\ref{sec:review_adet}. In what follows we concentrate on describing and analyzing the triggering mechanism for the communication from controller to sensors. We present first the kind of estimates of the norm of the state that the controller can obtain using the information it receives from the sensors. These estimates are then used to derive the triggering condition at the controller. At the end of the section we analyze the more practical aspects of the proposed implementation.

We start introducing the following assumption restricting the type of ISS controllers amenable to the strategy we propose in what follows:
\begin{assum}
\label{ass:growth}
The ISS closed-loop system~(\ref{eq:sys_iss}), satisfies the following properties:
\begin{equation} \lim_{s\to 0} \underline{\alpha}^{-1}\circ\overline{\alpha}(\underline\alpha^{-1}\circ\alpha_v^{-1}\circ\alpha_e(s)+2s)s^{-1}<\infty\end{equation}
\end{assum}
\begin{rem}
Note that this assumption, as well as  as Assumptions~\ref{ass:ISS} and~\ref{ass:lipschitz}, are automatically satisfied by linear systems with a linear state feedback controller and  the usual (ISS) Lyapunov function $V(x)=\sqrt{x^TPx}$.
\end{rem}

\subsection{Bounds of the system's state}
\label{ssec:bounds}

We start by establishing bounds of the norm of the state of the system computable at the controller. 
Denote by $\hat{\xi}$ the vector with entries:
$$\hat{\xi}_i(t)=\xi_i(t^i_{r_i}),\, t\in[t^i_{r_i}, t^i_{r_i+1}[\,\forall\,i=1,\ldots,n.$$ 
Thus, $\hat{\xi}$ can be seen as the value of the state vector that the controller is using to compute the input to the system.
The controller can compute then the upper bound: 
\begin{eqnarray}\notag
|\overline{\xi}|(t)&:=&|\hat{\xi}(t)|+\eta(t^c_{r_c})\geq |\hat{\xi}(t)-\varepsilon(t)|\\
&=&|\xi(t)|,\; \forall\, t\in[t^c_{r_c},t^c_{r_c+1}[,
\end{eqnarray}
which also satisfies the bound $|\overline{\xi}|(t)<|\xi(t)|+2\eta(t^c_{r_c})$.

\subsection{Triggering condition for the update of thresholds}

The following theorem proposes a condition to trigger the update of sensor thresholds guaranteeing UGAS of the closed-loop:

\begin{thm}
\label{thm:UGAS}
Consider the closed loop system~(\ref{eq:system}), (\ref{eq:controller}), (\ref{eq:eta_variable}), (\ref{eq:upd_rule_dyn}) with the threshold update rule~(\ref{eq:eta_upd_rule}). If Assumptions~\ref{ass:ISS}, \ref{ass:lipschitz} and \ref{ass:growth} hold, there exists $\rho<\infty$ sufficiently large, such that employing the sequence of threshold update times $\lbrace t^c_{r_c}\rbrace$, implicitly defined by:
\begin{eqnarray}
\label{eq:thres_trigger}
t^c_{r_c+1}&:=&\min\lbrace t=t^c_{r_c}+r\tau^c\,|\,r\in\N^+,\\ \notag 
&&|\overline{\xi}|(t)\leq\overline{\alpha}^{-1}\circ\underline{\alpha}(\rho\eta(t^c_{r_c}))\rbrace,
\end{eqnarray}
and setting:
\begin{eqnarray}
\label{eq:thres_0}
\eta(t^c_0)&\geq&\frac{\mu}{\rho}\underline\alpha^{-1}\circ V(\xi(t^c_0))
\end{eqnarray}
renders the closed-loop UGAS.
\end{thm}
\begin{pf}
We use Lemma~\ref{lem:UGAS} to show the desired result. The first itemized condition of the lemma is satisfied by the employment of the update rule~(\ref{eq:eta_upd_rule}) with a constant $\mu\in]0,1[$. Thus, we only need to show that the sequence $\lbrace t^c_{r_c}\rbrace$ is divergent and that the second itemized condition in the lemma also hold.  

First we show that starting from some time $t^c_{r_c}$ there always exists some time $T\geq t^c_{r_c}$ such that for all $t\geq T$ $|\overline{\xi}|(t)\leq\overline{\alpha}^{-1}\circ\underline{\alpha}(\rho\eta(t^c_{r_c}))$. Showing this guarantees that $\lbrace t^c_{r_c}\rbrace$ is a divergent sequence. 
From Assumption~\ref{ass:ISS} we know that for every $\epsilon>0$ there exists some $T\geq t^c_{r_c}$ such that $\underline\alpha(|\xi(t)|)\leq V(\xi(t))\leq\alpha_v^{-1}\circ\alpha_e(\eta(t^c_{r_c}))+\epsilon$ for every $t\geq T$. And therefore: $$|\overline\xi|(t)\leq \underline\alpha^{-1}(\alpha_v^{-1}\circ\alpha_e(\eta(t^c_{r_c}))+\epsilon)+2\eta(t^c_{r_c})$$ for every $t\geq T$. Thus, if there exists a $\rho>0$ such that
\begin{equation}
\underline{\alpha}^{-1}\circ\overline{\alpha}(\underline\alpha^{-1}(\alpha_v^{-1}\circ\alpha_e(s)+\epsilon)+2s)\leq\rho s
\end{equation}
holds for all $s\in]0,\eta(t_0)]$ and some $\epsilon>0$, a triggering event will eventually happen. Note that as $\epsilon$ can be made arbitrarily small (at the cost of possibly arbitrarily large, but finite, $T$), and due to the continuity of $\mathcal{K}_\infty$ functions, one can simply require that there exists a $\rho<\infty$ such that
\begin{equation}
\label{eq:rho_bound}
\underline{\alpha}^{-1}\circ\overline{\alpha}(\underline\alpha^{-1}\circ\alpha_v^{-1}\circ\alpha_e(s)+2s)<\rho s
\end{equation}
holds for all $s\in]0,\eta(t_0)]$. Finally, from Assumption~\ref{ass:growth} and Lemma~\ref{lem:growth_bound} one can conclude that such a $\rho<\infty$ exists.

The second condition of Lemma~\ref{lem:UGAS} is easier to prove. We start remarking that with $\rho$ so that~(\ref{eq:rho_bound}) holds, as $\mu<1$ and $\underline\alpha^{-1}\circ\overline\alpha(s)>s$ for all $s$, we also have: 
$$
\frac{\rho}{\mu}s>\underline\alpha^{-1}\circ\alpha_v^{-1}\circ\alpha_e(s)
$$
for all $s\in]0,\eta(t_0)]$.
Note next that as $V(\xi(t^c_{r_c}))\leq\overline\alpha(|\overline\xi|(t^c_{r_c}))$, from the triggering condition~(\ref{eq:thres_trigger}) and (\ref{eq:thres_0}), at all times $t^c_{r_c}$ the following holds: $\underline\alpha^{-1}\circ V(\xi(t^c_{r_c})\leq \frac{\rho}{\mu}\eta(t^c_{r_c})$. Therefore
\begin{eqnarray}
\kappa:=\frac{\rho}{\mu}\geq \frac{\underline{\alpha}^{-1}(\max\lbrace V(\xi(t^c_{r_c})),\alpha_v^{-1}\circ\alpha_e(\eta(t^c_{r_c})) \rbrace)}{\eta(t^c_{r_c})}
\end{eqnarray}
for all $t^c_{r_c}$, which concludes the proof.

\end{pf}
\begin{rem}
\label{rem:growth}
Finding controllers to satisfy Assumption~\ref{ass:growth} might be, in general, an arduous task. However, in practice one is generally only concerned with attaining practical stability and thus disregard this assumption. Then it is enough to guarantee~(\ref{eq:rho_bound}) for $s\in[\eta_{min},\eta(t^c_0)]$, $\eta_{min}>0$, which can always be satisfied given that: for every $\alpha\in\mathcal{K}_\infty$ there always exists $\kappa<\infty$ such that $\alpha(s)\leq \kappa s$ for all $s\in [\eta_m,\eta(t^c_0)]$.
\end{rem}

\begin{cor}
\label{cor:cor1}
If global bounds exists of the form: $$\overline\alpha(s)\leq\overline\kappa s,\, \overline\alpha(s)\geq\underline\kappa s,\,\alpha_v^{-1}\circ\alpha_e(s)\leq \kappa_{ve}s,$$ then any $$\rho>\frac{\overline\kappa}{\underline\kappa^2}(\kappa_{ve}+2\underline\kappa)$$ is sufficiently large to attain UGAS.
\end{cor}

\subsection{One-bit-communications}

We discuss now the implementation of asynchronous event-triggered controllers using one bit communications. The observation that allows such implementations is rather simple: to recover the value of a sensor after a threshold crossing, it is only necessary to know the previous value of the sensor and the sign of the error $\varepsilon_i$ when it crossed the threshold. In fact, if one assumes an initial round of synchronized measurements when the system is initialized, and these are transmitted completely (with enough bits), then the necessary piece of information to recover the state can be represented by a single bit indicating the sign of the error. More specifically, the actual sensor measurement can be recovered recursively as:
\begin{eqnarray}
\label{eq:hat_xi}\notag
 \hat{\xi}_i(t^i_{r_i})&=&\hat{\xi}_i(t^i_{r_i-1})+(2d_i(t^i_{r_i})-1)\sqrt{\eta_i(t^c_{r_c})},\, t^i_{r_i}\in[t^c_{r_c}, t^c_{r_c+1}[\\
 \hat{\xi}_i(t^i_0)&=&\xi(t_0),
\end{eqnarray}
where $d_i(t^i_{r_i})\in\lbrace{0,1\rbrace}$ denotes the negative ($0$) or positive ($1$) sign of the error $\varepsilon_i$ when the triggering event was released. Remember that the controller can actually keep track of the values $\eta_i$ employed by each of the sensors.
Similarly, the messages from the controller to the sensors commanding a reduction of the thresholds can be indicated with a single bit, which would also leave the possibility open for commanding the increment of thresholds or a synchronized measurement e.g. whenever the controller detects the system suffered an impulsive disturbance.

\subsection{Global inter-transmission bounds}

In order to make the technique proposed in this section practical, we still need to guarantee that there exists a minimum time between sensor transmissions. Theorem~\ref{thm:UGAS}, by means of Lemma~\ref{lem:UGAS}, guarantees that between two threshold updates there exists a minimum time between sensor transmissions of the same sensor, \ie:
\begin{equation}
\label{eq:bas_t_bound}
t^i_{r_i}-t^i_{r_i-1}\geq \tau_i^*:=\frac{L_{f_i}^{-1}\theta_i}{1+\kappa},\; \forall\,t^i_{r_i}, t^i_{r_i-1}\in [t^c_{r_c},t^c_{r_c+1}[,
\end{equation}
with $\kappa:=\frac{\rho}{\mu}$ (see proof of Theorem~\ref{thm:UGAS}).

However, it could happen that at times $t^c_{r_c}$ some sensors, by reducing their local threshold, automatically violate their triggering condition. Further, this could lead to times between transmissions of the same sensor violating the bounds from Lemma~\ref{lem:bounds}. Fortunately, there is a simple solution guaranteeing new lower bounds of the time between transmissions.

\begin{prop}
Let $\mu\in ]\sqrt{0.5},1[$ and the local thresholds be modified as:
\begin{equation}
\label{eq:loc_thresh_upd}
\eta_i(t):=\theta_i^2\eta^2(t-\tilde{\tau}^*_i),\,|\theta|=1,
\end{equation}
{\it i.e.} update the local thresholds only $\tilde{\tau}^*_i$ units of time after the controller detects that the thresholds can be reduced.  

The closed loop system~(\ref{eq:system}), (\ref{eq:controller}), (\ref{eq:eta_variable}), (\ref{eq:loc_thresh_upd}) with the threshold update rule~(\ref{eq:eta_upd_rule}), (\ref{eq:thres_trigger}) is UGAS, and there exists a minimum time between events at each sensor given by:
\begin{small}
\begin{equation}
\label{eq:bounds_corrected}
\tilde{\tau}_i^*:=\left(1-\sqrt{\frac{1-\mu^2}{\mu^2}} \right)\tau_i^*.
\end{equation}
\end{small}
\end{prop}

\begin{pf}
The proof for UGAS is identical to that provided for Theorem~\ref{thm:UGAS}, as the only change is that the actual sequence of times at which the thresholds are effectively updated is shifted by $\tilde{\tau}_i^*$. Also, as stated earlier, all updates happening between two threshold updates respect the bound~(\ref{eq:bas_t_bound}). Thus, we only need to study the inter transmission times due to triggering at the times $t^c_{r_c}$.
Consider some threshold update time $t^c_{r_c}$ and let sensor $i$ have an error signal $\varepsilon_i^2(t^c_{r_c})>\eta_i(t^c_{r_c})$. One can always bound how much $\varepsilon_i^2$ can be above the new threshold as:
$$
\varepsilon_i^2(t^c_{r_c})-\eta_i(t^c_{r_c})\leq \eta_i(t^c_{r_c-1})-\eta_i(t^c_{r_c}) =(\frac{1}{\mu^2}-1)\eta_i(t^c_{r_c}).
$$ 
If the reduction of the threshold triggers a sensor transmission, \ie $t^i_{r_i}=t^c_{r_c}$ because $\varepsilon_i^2(t^c_{r_c})>\eta_i(t^c_{r_c})$, then
the controller is updated with a value:
\begin{equation}
\label{eq:sign}
\hat{\xi}_i(t^i_{r_i})=\hat{\xi}_i(t^i_{r_i-1})+\texttt{sign}(\varepsilon_i(t^i_{r_i}))\sqrt{\eta_i(t^c_{r_c})}.
\end{equation}
This is equivalent to reseting the local error at the triggering sensor, after such a controller update, to a non-zero value:
$$
\varepsilon_i^2(t^c_{r_c})\in [0,\; (\frac{1}{\mu^2}-1)\eta_i(t^c_{r_c})[.
$$
Thus, reasoning as in the proof of Lemma~\ref{lem:bounds}, but now computing the time it takes $|\varepsilon_i|$ to go from a value of $\sqrt{(\frac{1}{\mu^2}-1)}\sqrt{\eta_i(t^c_{r_c})}$ to $\sqrt{\eta_i(t^c_{r_c})}$, one can show that \mbox{$t^i_{r_i+1}-t^i_{r_i}\geq \tilde\tau_i^*$}, whenever $t^i_{r_i}=t^c_{r_c}$. Furthermore, the use of the local thresholds~(\ref{eq:loc_thresh_upd}) guarantees that also $t^i_{r_i}-t^i_{r_i-1}\geq \tilde\tau_i^*$, whenever $t^i_{r_i-1}\in ]t^c_{r_c-1}, t^c_{r_c}[$ and $t^i_{r_i}\in ]t^c_{r_c}, t^c_{r_c+1}[$.
Finally, note that $\mu>\sqrt{0.5}$ is required to provide a lower bound $\tilde\tau_i^*$ that is positive.
\end{pf}

\subsection{Other practical issues}
\label{ssec:practical}

\subsubsection{Delays}
Many effects of a real practical implementation can be abstracted in the form of a delay in the proposed event-triggered implementation.
We illustrate this with a specific example: In our implementations controller updates can take place arbitrarily close to each other. This is so because while one sensor cannot trigger updates arbitrarily close to each other, the combination of all sensors can potentially force that to happen. This makes the proposed techniques more suitable for systems with controller(s) and actuators co-located. As we did in~\cite{MazoCaoCDC2011}, we suggest the use of a periodic subjacent scheme for the update of the controller. The effect of such a scheme is the introduction of an artificial delay in the closed-loop system.

The kind of delays we consider are those between the event-generation at the sensor side and its effect taking place in the control inputs applied to the system. Essentially, what most event-triggered techniques do is control the magnitude of the virtual error introduced by sampling in a digital implementation. If the magnitude of this error signal is successfully kept within certain margins, the controller implementation is stable. This error signal that one must control is defined at the plant side. Therefore, when delays are present, while the sensors send new measurements trying to keep $|\varepsilon_i(t)|=|\xi_i(t^i_{r_i})-\xi_i(t)|\leq \eta_i$, what actually matters is the value of the error at the plant-side $\hat{\varepsilon}_i(t)$, defined as:
\begin{eqnarray}
\hat{\varepsilon}_i(t)&=&\xi(t^i_{r_i-1})-\xi(t)\qquad t\in [t^i_{r_i}, t^i_{r_i}+\Delta\tau^i_{r_i}[\\
\hat{\varepsilon}_i(t)&=&\varepsilon_i(t)\qquad\qquad\qquad\; t\in [t^i_{r_i}+\Delta\tau^i_{r_i}, t^i_{r_i+1}[,
\end{eqnarray}
where $\Delta\tau^i_{r_i}$ denotes the delay between the time $t^i_{r_i}$ at which a measurement is transmitted and the time $t^i_{r_i}+\Delta\tau^i_{r_i}$ at which the controller is updated with that new measurement. Thus the actual objective to attain UGAS or UGPS is to keep $|\hat\varepsilon(t)|$ below the threshold $\eta$.
From the analysis in the proof of Lemma~\ref{lem:bounds} we know that the maximum speed of the error signal is always kept below $L_{f_i}(\kappa+1)\eta$
and thus, given a maximum delay of $\Delta\tau^i_{r_i}\leq\Delta\tau<(\kappa+1)^{-1}L_{f_i}^{-1}$ for all $i$ and $r_i$, reducing the local thresholds as:
$$\theta_i\bar{\eta}=\theta_i\eta\left(1-L_{f_i}(\kappa+1)\Delta\tau\right)$$ and keeping $|\varepsilon(t)|\leq\bar\eta$, guarantees that the error at the plant side stays below the desired value $|\hat{\varepsilon}_i(t)|\leq\eta$. The more conservative our estimates of $\kappa$ and $L_{f_i}$ are, the smaller the tolerable delays will be.

\subsubsection{Disturbances}

Another important question is how to deal with the presence of disturbances or sensing/actuation noise.
Whenever a persistent disturbance is present, as could be sensing noise, in general one can only achieve practical stability to a region dependent of the power of the disturbance. 
Consider for the discussion the closed loop dynamics $\dot{\xi}=f(\xi,k(\xi+\varepsilon),\delta)$, where $\delta$ is a disturbance signal. Now, further assume the closed-loop system is ISS also with respect to the disturbance, \ie there exists an ISS Lyapunov function $V$ such that $\dot{V}(x,e,d)\leq -\alpha_v\circ V(x)+\alpha_e(|e|)+\alpha(|d|)$. In order to obtain minimum inter-transmission times greater than zero, it is obvious that now we cannot let $\eta$ go to zero. Instead, one would have to select a minimum $\eta_m$ guaranteeing a region to which the system converges small enough for the application at hand. If the disturbance has bounded magnitude, \ie $\Vert\delta\Vert_\infty\leq n_d\eta_m$ with some $n_d\in\R^+$, one can easily show that the minimum inter-transmission times guaranteed by Lemma~\ref{lem:bounds} can be recomputed to be: $$\tau_i^*= \frac{L_{f_i}^{-1}\theta_i}{(\kappa+1)+n_d}.$$
Furthermore, imposing such a minimum threshold carefully one can also prevent that measurement noise forces a triggering on its own (which could lead to a wrong estimation of the state at the controller).
If one would like to design a system resilient also to impulsive disturbances, rather than employing the analysis just proposed, it might be a better idea to enable a mechanism in the implementation to force a synchronous update (essentially a reinitialization of the system) to compensate for such impulses. A more detailed analysis of these solutions is left for future work.

\subsubsection{Performance guarantees}

Finally, we would like to make some remarks with respect to the performance guarantees. While it may be possible to provide explicit performance guarantees of the implementation (see the proof of Theorem~\ref{thm:UGAS}), it would take a very complicated form rendering it arguably practical. Instead of providing such performance guarantees, let us just give some intuition on how the different design parameters affect the convergence speed of the implementation and the triggering of events at sensors and controller. There are three design parameters for the proposed implementation: $\rho$, $\mu$ and $\tau^c$. Reducing any of these parameters in general should yield a faster convergence of the system. However, while reducing $\mu$ in exchange may lead to more frequent inter-transmission times from the sensors, reducing $\rho$ may increase the frequency of threshold update requests sent from the controller to the sensors.

\section{Examples}
\label{sec:example}

\subsection{Example 1}

\begin{figure}[ht]
 \centering
 \includegraphics[width=\hsize]{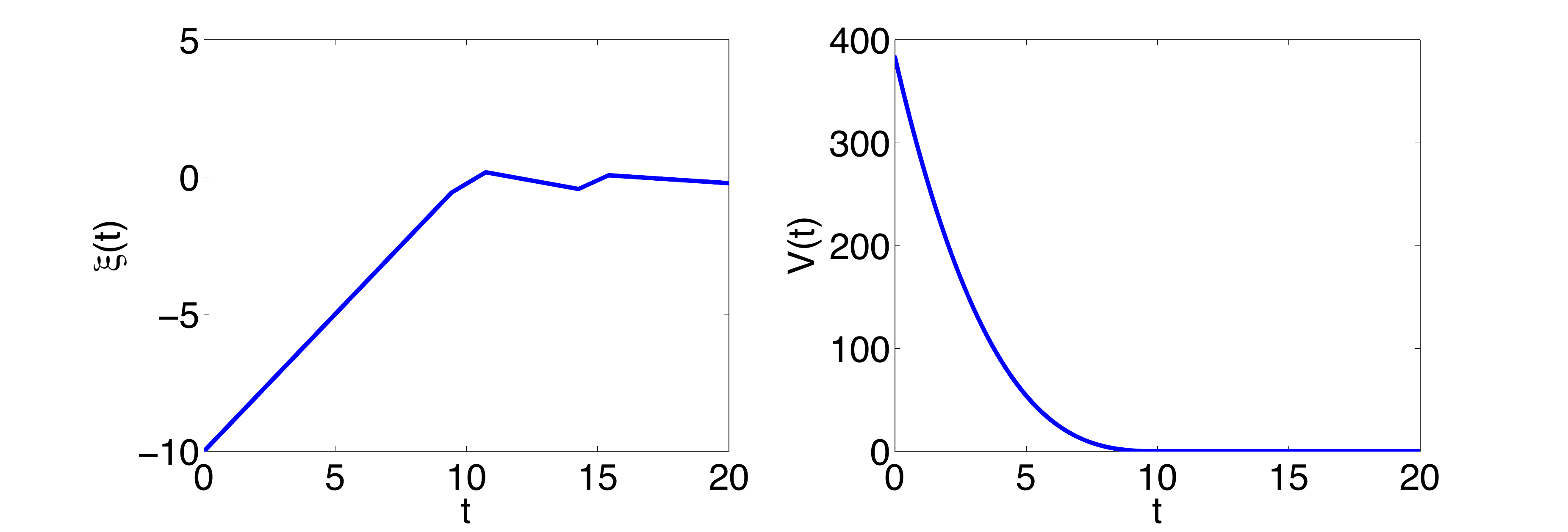}\\
 \includegraphics[width=\hsize]{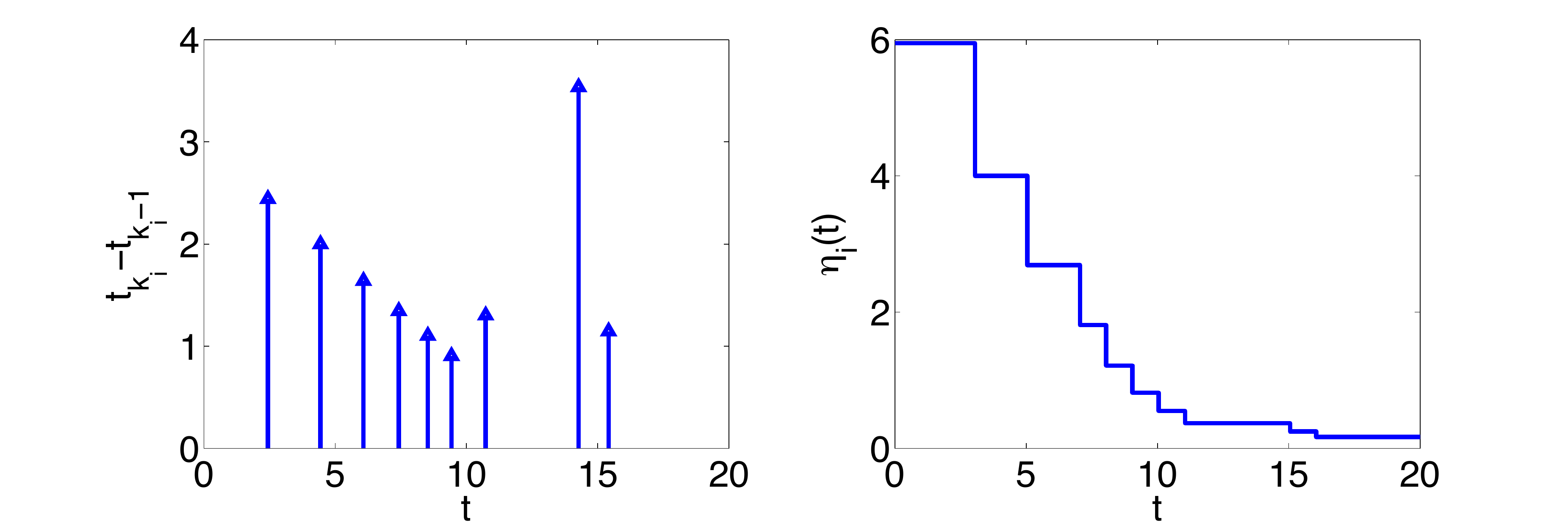}
 \caption{State trajectory, Lyapunov function evolution, events generated at the sensor and evolution of the threshold.}
 \label{fig:fig_sat1}
 \end{figure}

Let us employ for simplicity a scalar example to illustrate the role of Assumption~\ref{ass:growth} in the design of the proposed implementation. 
Consider the system:
\begin{equation}
\dot{\xi}(t)=\texttt{sat}(\upsilon(t))
\end{equation}
with a controller affected by measurement errors: $\upsilon(t)=-\xi(t)-\varepsilon(t)$, where the function $\texttt{sat}(s)$ is the saturation function, saturating at values $|s|>1$. 
In~\cite{sontag1995state} it is shown that with an ISS Lyapunov function $V(x)=\frac{|x|^3}{3}+\frac{|x|^2}{2}$, the following dissipation inequality holds:
\begin{equation}
\dot{V}(x,e)\leq -\frac{x^2}{2}+2e^2.
\end{equation}
In this particular example one can take $\overline\alpha(s)=\underline\alpha(s)=s^3/3+s^2/2$. Furthermore, using $\alpha_v(s)=\alpha_x\circ\overline\alpha^{-1}(s)$, with $\alpha_x(s)=s^2/2$, and $\alpha_e(s)=2s^2$:
\begin{equation}
 \underline{\alpha}^{-1}\circ\overline{\alpha}(\underline\alpha^{-1}\circ\alpha_v^{-1}\circ\alpha_e(s)+2s) = \alpha_x^{-1}\circ\alpha_e(s)+2s=4s,
\end{equation}
and we can conclude that any $\rho>4$ will guarantee asymptotic stability to the origin. 
Furthermore, selecting \eg $\rho=4.1$,  $\mu=0.82$ and $\tau^c=1$, results in a minimum time between sensor transmissions of $0.04$, after accounting for the effect of an aggregated (communication/actuation) delay of $0.002$. We show in Figure~\ref{fig:fig_sat1} the result of a simulation, where $\xi(0)=-10$ and $\varepsilon(0)=0$, in which it can be appreciated how the system is stabilized while the inter-transmission times respect the lower-bound.

\subsection{Example 2}

\begin{figure}[ht]
 \centering
 \includegraphics[width=\hsize]{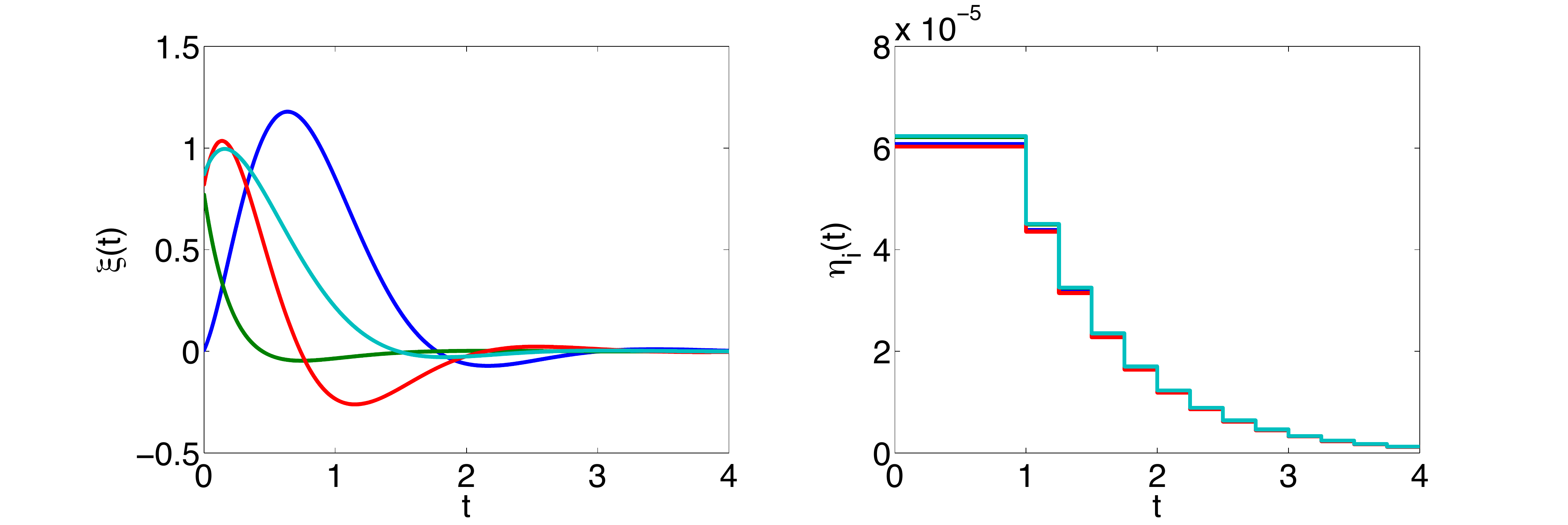}\\
 \includegraphics[width=\hsize]{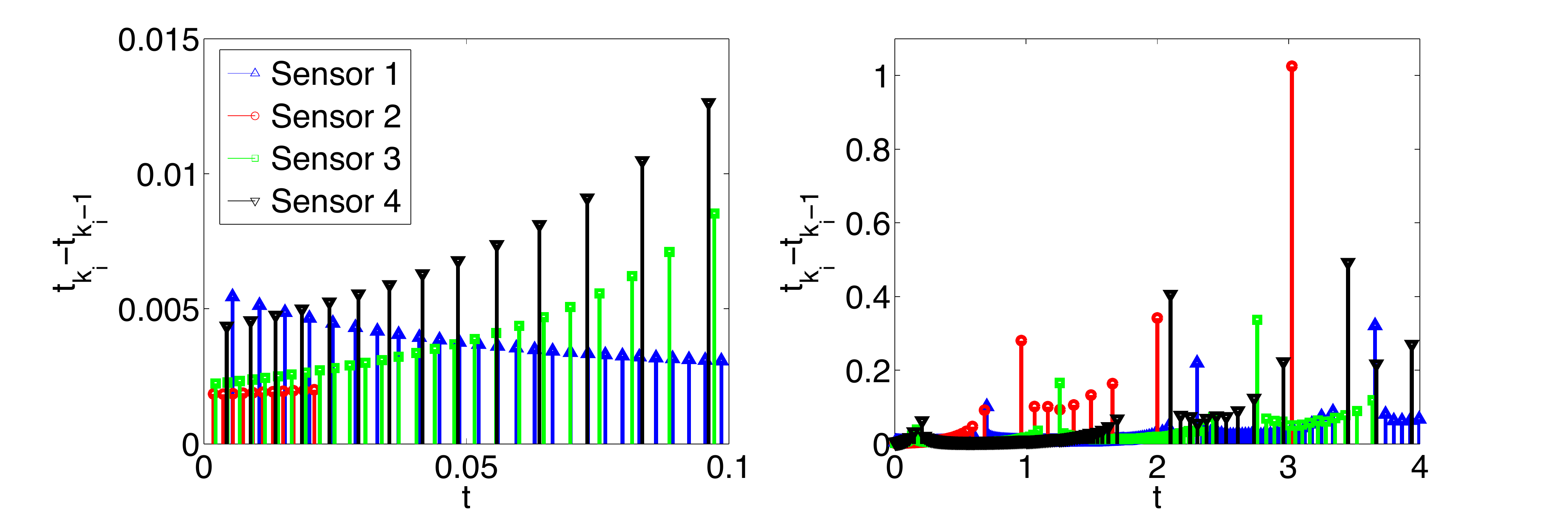}
 \caption{State trajectory, evolution of the thresholds, and events generated at the sensors.}
 \label{fig:fig_qlin1}
 \end{figure}

In this next example we illustrate the asynchrony of measurements. Consider a nonlinear system of the form:
\begin{equation}
\label{eq:sys_ex}
\dot{\xi}(t)=A\xi +B(f(\xi(t))+\upsilon(t))
\end{equation}
where $f$ is a nonlinear locally Lipschitz function.
One can stabilize these kind of systems employing a feedback of the form: $\upsilon(t)=-f(\xi(t))-K\xi(t)$ with $K$ such that $A_c=A-BK$ is Hurwitz.
Consider now the feedback affected by measurement errors: $$\upsilon(t)=-f(\xi(t)+\varepsilon(t))-K(\xi(t)+\varepsilon(t)),$$
and let $V(x)=\sqrt(x^TPx)$, where $PA_c+A_c^TP=-I$, be the candidate ISS-Lyapunov function for the system.
This function is indeed an ISS-Lyapunov function as it satisfied the following dissipation inequality:
\begin{eqnarray*}
\dot{V}(x,e)\leq-\frac{1}{2\sqrt{\lambda_M(P)}}|x| +\left(|\sqrt{P}BK|+|\sqrt{P}|L_f\right)|e|,
\end{eqnarray*}
and where $L_f$ is the Lipschitz constant of $f$ in the compact where we assume the system will evolve (e.g. computed from the compact set to which the initial condition of the system is assumed to belong).
Furthermore, note that $\sqrt{\lambda_m(P)}|x|\leq V(x) \leq \sqrt{\lambda_M(P)}|x|$.
Thus, the $\mathcal{K}_\infty$ functions $\underline{\alpha}$, $\overline{\alpha}$, $\alpha_e$ and $\alpha_v$ can be assumed to be linear functions in the compact of interest. This guarantees that Assumption~\ref{ass:growth} is satisfied. Furthermore, by Collorary~\ref{cor:cor1}, a $\rho$ that would guarantee Asymptotic stability is:
\begin{scriptsize}
\begin{equation}
\rho>2\frac{\lambda_M(P)}{\lambda_m(P)}\left( \sqrt{\lambda_M(P)}\left(|\sqrt{P}BK|+|\sqrt{P}|L_f\right)+\sqrt{\frac{\lambda_M(P)}{\lambda_m(P)}}\right).
\end{equation}
\end{scriptsize}
We provide simulations now for a particular system of this form, with:
\begin{scriptsize}
\begin{eqnarray}
A=\left[\begin{array}{cccc} 1.5 & 0 & 7 &  -5\\
   -0.5 & -4 & 0 & 0.5\\
    1  & 4 & -6 & 6\\
    0  & 4 & 1 & -2
\end{array} \right],
\;
B=\left[\begin{array}{cc} 0 & 0\\ 5 & 0\\ 1 & -3\\1 & 0
\end{array} \right],
\\
f(x)=\left[\begin{array}{c}
x_2^2\\\sin(x_3)
\end{array} \right],\;
K=\left[\begin{array}{cccc} 0.1 & -0.2 & 0 & -0.2\\ 1.5 & -0.2 & 0 & 0
\end{array} \right],
\end{eqnarray}
\end{scriptsize}
In Figure~\ref{fig:fig_qlin1} it is shown the result of a simulation with the following design parameters: $\mu=0.85$, $\rho=253$, $\tau^c=0.25$, initialized with $\xi(0)=[0\,0.8\,0.7\,0.75]^T$ and $\varepsilon(0)=0$. The function $f(x)$ is only locally Lipschitz, which in order to guarantee asymptotic stability for initial conditions $|\xi(t_0)|\leq 2$ imposed a $\rho>230$. With this design, and assuming there is a maximum delay introduced of $\Delta\tau=2\cdot10^{-5}$, results in a minimum time between transmissions of the same sensor of $5\cdot10^{-5}$. The simulation shows that these bounds are conservative, as the minimum observed inter transmission time was $0.0018$. Furthermore, the average inter transmission time in the simulated time was one order of magnitude larger.

\section{Discussion}
\label{sec:discuss}

We have shown how asymptotic stability can be attained with fully decentralized event-triggered implementations. At this point, it is important to make a distinction between having decentralized event-triggered conditions, and having a decentralized controller. We have concentrated on the process of deciding when to transmit measurements between sensors, controllers and actuators, making it truly decentralized. It would also be interesting to see how the present approach combines with the case of decentralized/distributed controllers~\cite{Wang11}.

Some discussion needs to be provided with respect to the assumptions put in place. The main of the assumptions being the requirement of having a controller rendering the system ISS on compacts with respect to measurement errors. It is important to remark that this property does not need to be globally satisfied but only in the compact of interest for the system operation, which relaxes drastically the requirement and makes it also easier to satisfy~\cite{Freeman95}. Furthermore, in practice one usually is only concerned with practical stability, which as stated in Remark~\ref{rem:growth} eliminates the need for Assumption~\ref{ass:growth}.

As the simulated example illustrated, the guarantees that the analysis provides are highly conservative. The main reason for this conservatism can be found in the use of Lipschitz constants for the bounding of the speed of the system on a compact. While the assumptions put in place provide with a clean theory and easy to follow results, it would be desirable to study computational methods capable of reducing this conservatism gap. This idea is not new, and has already been applied to several techniques in the case of linear systems employing LMI solvers~\cite{DonkersLMI,HetelLMI,Donkers12}.

Following these comments we leave as follow-up work the design of computational methods to relax the conservative bounds obtained, and also to help in the design itself of the implementation, \ie the selection of adequate $\rho$ and $\mu$.This approach also opens the door to perform analysis to optimize controller design and implementation for data-rate transmissions. Also, in the more practical side of the future work, it is desirable the study and design of protocols for wireless communications exploiting the benefits of the proposed techniques. On the more theoretical side, many of the practical effects briefly discussed in Section~\ref{ssec:practical} can, and probably should, be studied in more detail. In particular, and connected to the computational approaches, an important issue to be solved is the provision of guarantees of performance, including the response to disturbances. Finally, it would be interesting to study controller designs for useful classes of non-linear systems that satisfy the current assumptions, or relaxed versions of them.

\bibliographystyle{plain}
\bibliography{cdc11} 

\begin{thebibliography}{10}

\bibitem{Sontag:2008p2453}
A.~Agrachev, A.~S. Morse, E.~Sontag, H.~Sussmann, and V.~Utkin.
\newblock Input to state stability: Basic concepts and results.
\newblock In {\em Nonlinear and Optimal Control Theory}, volume 1932 of {\em
  Lecture Notes in Mathematics}, pages 163--220. Springer Berlin / Heidelberg,
  2008.

\bibitem{AntaTAC}
A.~Anta and P.~Tabuada.
\newblock To sample or not to sample: Self-triggered control for nonlinear
  systems.
\newblock {\em IEEE Transactions on Automatic Control}, 55:2030--2042, Sep.
  2010.

\bibitem{AraujoDCOSS}
J.~Araujo, A.~Anta, M.~Mazo~Jr., J.~Faria, A.~Hernandez, P.~Tabuada, and K.H.
  Johansson.
\newblock Self-triggered control over wireless sensor and actuator networks.
\newblock In {\em Distributed Computing in Sensor Systems and Workshops
  (DCOSS), 2011 International Conference on}, pages 1 --9, june 2011.

\bibitem{Astrom}
K.J. {\AA}str\"{o}m and B.M. Bernhardsson.
\newblock Comparison of {Riemann} and {Lebesgue} sampling for first order
  stochastic systems.
\newblock In {\em Proceedings of the 41st IEEE Conference on Decision and
  Control}, volume~2, pages 2011 -- 2016, Dec. 2002.

\bibitem{Cervin08}
A.~Cervin and T.~Henningsson.
\newblock Scheduling of event-triggered controllers on a shared network.
\newblock In {\em Proceedings of the 47th IEEE Conference on Decision and
  Control}, pages 3601 --3606, Dec. 2008.

\bibitem{Dimarogonas12}
D.V. Dimarogonas, E.~Frazzoli, and K.H. Johansson.
\newblock Distributed event-triggered control for multi-agent systems.
\newblock {\em Automatic Control, IEEE Transactions on}, 57(5):1291--1297,
  2012.

\bibitem{DonkersLMI}
M.~Donkers, W.~Heemels, N.~van~de Wouw, and L.~Hetel.
\newblock {Stability analysis of networked control systems using a switched
  linear systems approach}.
\newblock {\em IEEE Transactions on Automatic Control}, 56(9):2101--2115, Sept.
  2011.

\bibitem{Donkers12}
M.C.F. Donkers and W.~Heemels.
\newblock Output-based event-triggered control with guaranteed l∞-gain and
  improved and decentralized event-triggering.
\newblock {\em IEEE Transactions on Automatic Control}, 57(6):1362--1376, 2012.

\bibitem{Freeman95}
R.~Freeman.
\newblock Global internal stabilizability does not imply global external
  stabilizability for small sensor disturbances.
\newblock {\em IEEE Transactions on Automatic Control}, 40(12):2119 --2122, dec
  1995.

\bibitem{Heemels08}
W.P.M.H. Heemels, J.H. Sandee, and P.P.J. van~den Bosch.
\newblock {Analysis of event-driven controllers for linear systems}.
\newblock {\em International Journal of Control}, 81(4):571--590, 2008.

\bibitem{HetelLMI}
L.~Hetel, A.~Kruszewski, W.~Perruquetti, and J.~Richard.
\newblock Discrete and intersample analysis of systems with aperiodic sampling.
\newblock {\em IEEE Transactions on Automatic Control}, 56(7):1696 --1701, july
  2011.

\bibitem{DynamicQuantizers}
D.~Liberzon.
\newblock Hybrid feedback stabilization of systems with quantized signals.
\newblock {\em Automatica}, 39(9):1543--1554, 2003.

\bibitem{Lunze10}
J.~Lunze and D.~Lehmann.
\newblock A state-feedback approach to event-based control.
\newblock {\em Automatica}, 46(1):211 -- 215, 2010.

\bibitem{MazoCaoCDC2011}
M.~Mazo~Jr. and M.~Cao.
\newblock Decentralized event-triggered control with asynchronous updates.
\newblock In {\em Proceedings of the 50th Conference on Decision and Control},
  pages 2547--2552, 2011.

\bibitem{MazoCaoADHS2012}
M.~Mazo~Jr. and M.~Cao.
\newblock Decentralized event-triggered control with one bit communications.
\newblock In {\em 4th IFAC Conference on analysis and design of hybrid systems
  (ADHS'12), Eindhoven, The Netherlands}, pages 52--57, 2012.

\bibitem{MazoTACWSAN}
M.~Mazo~Jr. and P.~Tabuada.
\newblock Decentralized event-triggered control over wireless sensor/actuator
  networks.
\newblock {\em IEEE Transactions on Automatic Control, Special issue on
  Wireless Sensor Actuator Networks}, 56(10):2456--2461, Oct. 2011.

\bibitem{Molin10}
A.~Molin and S.~Hirche.
\newblock Optimal event-triggered control under costly observations.
\newblock In {\em Proceedings of the 19th International Symposium on
  Mathematical Theory of Networks and Systems}, pages 2203--2208, 2010.

\bibitem{DeadbandTilbury}
P.G. Otanez, J.R. Moyne, and D.M. Tilbury.
\newblock Using deadbands to reduce communication in networked control systems.
\newblock In {\em American Control Conference, 2002. Proceedings of the 2002},
  volume~4, pages 3015 -- 3020 vol.4, 2002.

\bibitem{Rabi08}
M.~Rabi and K.~H. Johansson.
\newblock Event-triggered strategies for industrial control over wireless
  networks.
\newblock In {\em Proceedings of the 4th Annual International Conference on
  Wireless Internet}, WICON '08, pages 34:1--34:7, ICST, Brussels, Belgium,
  Belgium, 2008. ICST (Institute for Computer Sciences, Social-Informatics and
  Telecommunications Engineering).

\bibitem{sontag1995state}
Eduardo Sontag.
\newblock State-space and i/o stability for nonlinear systems.
\newblock {\em Feedback control, nonlinear systems, and complexity}, pages
  215--235, 1995.

\bibitem{Tabuada}
P.~Tabuada.
\newblock Event-triggered real-time scheduling of stabilizing control tasks.
\newblock {\em IEEE Transactions on Automatic Control}, 52(9):1680 --1685,
  Sept. 2007.

\bibitem{Wang11}
X.~Wang and M.D. Lemmon.
\newblock Event-triggering in distributed networked control systems.
\newblock {\em IEEE Transactions on Automatic Control}, 56(3):586--601, 2011.

\bibitem{WeimerADHS12}
J.~Weimer, J.~Ara\'ujo, and K.H. Johansson.
\newblock Distributed event-triggered estimation in networked systems.
\newblock In {\em 4th IFAC Conference on analysis and design of hybrid systems
  (ADHS'12), Eindhoven, The Netherlands}, pages 178--185, 2012.

\bibitem{Estrin02}
Wei Ye, J.~Heidemann, and D.~Estrin.
\newblock An energy-efficient mac protocol for wireless sensor networks.
\newblock In {\em INFOCOM 2002. Twenty-First Annual Joint Conference of the
  IEEE Computer and Communications Societies. Proceedings. IEEE}, volume~3,
  pages 1567 -- 1576 vol.3, 2002.

\end{thebibliography}

\end{document}